\documentclass{article}
\usepackage[utf8]{inputenc}

\newtheorem{lemma}{Lemma}[section]

\newtheorem{theorem}[lemma]{Theorem}
\newtheorem{proposition}[lemma]{Proposition}
\newtheorem{corollary}[lemma]{Corollary}
\newtheorem{example}[lemma]{Example}

\newtheorem{note}[lemma]{Remark}

\title{Arities and aritizabilities of group, monoid and groupoid theories\footnote{The work
of the second author was carried out in the framework of the State
Contract of the Sobolev Institute of Mathematics, Project
No.~FWNF-2022-0012, and of Committee of Science in Education and
Science Ministry of the Republic of Kazakhstan (Grant No.
AP08855544).}}
\author{In.I. Pavlyuk, S.V. Sudoplatov}
\date{}

\begin{document}

\maketitle
\begin{abstract}
We study applications of a general approach for arities and
arizabilities of theories to group and monoid theories. It is
proved that a theory of a group $G$ is aritizable if and only if
$G$ is finite. It is shown that this criterion does not hold for
monoids/groupoids: there is an infinite monoid/groupoid having a
binary theory.
\end{abstract}


\bigskip
We continue to study arities of theories and of their expansions
\cite{aafot, ananat}. An arity is a basic characteristic of
complexity of a theory allowing to reduce all formulae to ones
with boundedly many free variables. In the present paper a general
approach for arities and arizabilities of theories \cite{aafot,
ananat} is applied to group theories \cite{Poizat}. It is proved
that a theory of a group $G$ is aritizable if and only if $G$ is
finite. It is shown that this criterion does not hold for
monoids/groupoids: there is an infinite monoid/groupoid having a
binary theory.

\section{Preliminaries}

Recall a series of notions related to arities and aritizabilities
of theories.

{\bf Definition} \cite{ZPS}. A theory $T$ is said to be
\emph{$\Delta$-based}\index{Theory!$\Delta$-based}, where $\Delta$
is some set of~formulae without parameters, if any formula of $T$
is equivalent in $T$ to a~Boolean combination of formulae
in~$\Delta$.

For $\Delta$-based theories $T$, it is also said that $T$ has {\em
quantifier elimination}\index{Elimination of quantifiers} or {\em
quantifier reduction}\index{Reduction of quantifiers} up to
$\Delta$.

\medskip
{\bf Definition} \cite{ZPS, CCMCT}. {\rm Let $\Delta$ be a set of
formulae of a theory $T$, and $p(\bar{x})$ a type of $T$ lying in
$S(T)$. The type $p(\bar{x})$ is said to be
\emph{$\Delta$-based}\index{Type!$\Delta$-based} if $p(\bar{x})$
is isolated by a set of formulas $\varphi^\delta\in p$, where
$\varphi\in\Delta$, $\delta\in\{0,1\}$.}

\medskip
The following lemma, being a corollary of Compactness Theorem,
noticed in \cite{ZPS}.

\medskip
\begin{lemma}\label{lem01} A  theory  $T$  is
$\Delta$-based  if and only if,  for  any  tuple~$\bar{a}$ of any
{\rm (}some{\rm )} weakly saturated model of $T$, the type ${\rm
tp}(\bar{a})$ is $\Delta$-based.
\end{lemma}

\medskip
{\bf Definition} \cite{aafot}. An elementary theory $T$ is called
{\em unary}, or {\em $1$-ary}, if any $T$-formula
$\varphi(\overline{x})$ is $T$-equivalent to a Boolean combination
of $T$-formulas, each of which is of one free variable, and of
formulas of form $x\approx y$.

For a natural number $n\geq 1$, a formula $\varphi(\overline{x})$
of a theory $T$ is called {\em $n$-ary}, or an {\em $n$-formula},
if $\varphi(\overline{x})$ is $T$-equivalent to a Boolean
combination of $T$-formulas, each of which is of $n$ free
variables.

For a natural number $n\geq 2$, an elementary theory $T$ is called
{\em $n$-ary}, or an {\em $n$-theory}, if any $T$-formula
$\varphi(\overline{x})$ is $n$-ary.

A theory $T$ is called {\em binary} if $T$ is $2$-ary, it is
called {\em ternary} if $T$ is $3$-ary, etc.

We will admit the case $n=0$ for $n$-formulae
$\varphi(\overline{x})$. In such a case $\varphi(\overline{x})$ is
just $T$-equivalent to a sentence
$\forall\overline{x}\varphi(\overline{x})$.

If $T$ is a theory such that $T$ is $n$-ary and not $(n-1)$-ary
then the value $n$ is called the arity of $T$ and it is denoted by
${\rm ar}(T)$. If $T$ does not have any arity we put ${\rm
ar}(T)=\infty$.

Similarly, for a formula $\varphi$ of a theory $T$ we denote by
${\rm ar}_T(\varphi)$ the natural value $n$ if $\varphi$ is
$n$-ary and not $(n-1)$-ary. If $\varphi$ does not any arity we
put ${\rm ar}_T(\varphi)=\infty$. If a theory $T$ is fixed we
write ${\rm ar}(\varphi)$ instead of ${\rm ar}_T(\varphi)$.

\begin{note} {\rm \cite{aafot} For the description of definable sets for models
$\mathcal{M}$ of $n$-theories it suffices describe links between
definable sets $A$ and $B$ for $n$-formulas
$\varphi(\overline{x})$ and $\psi(\overline{y})$, respectively,
and definable sets $C$ and $D$ for
$\varphi(\overline{x})\wedge\psi(\overline{y})$ and
$\varphi(\overline{x})\vee\psi(\overline{y})$, respectively.

If $\overline{x}=\overline{y}$ then $C=A\cap B$ and $D=A\cup B$,
i.e., conjunctions and disjunctions work as set-theoretic
intersections and unions.

If $\overline{x}$ and $\overline{y}$ are disjoint then $C=A\times
B$ and
$D=(A+B)_{\mathcal{M}}\rightleftharpoons\{\langle\overline{a},\overline{b}\rangle\mid\overline{a}\in
A\mbox{ and }\overline{b}\in M, \mbox{ or } \overline{a}\in
M\mbox{ and }\overline{b}\in B\}$, i.e., $C$ is the Cartesian
product of $A$ and $B$, and $D$ is the ({\em generalized}) {\em
Cartesian sum} of $A$ and $B$ in the model $\mathcal{M}$.

If $\overline{x}\ne\overline{y}$, and $\overline{x}$ and
$\overline{y}$ have common variables, then $C$ and $D$ are
represented as a {\em mixed product} and a {\em mixed sum},
respectively, working partially as intersection and union, for
common variables, and partially as Cartesian product and Cartesian
sum, for disjoint variables.

If $\overline{x}$ and $\overline{y}$ consist of pairwise disjoint
variables and $\overline{x}\subseteq\overline{y}$,
$\overline{x}\ne\overline{y}$, then for any formula
$\varphi(\overline{x})$ the set of solution of the formula
$\varphi(\overline{x})\wedge(\overline{y}\approx\overline{y})$ in
$\mathcal{M}$ is called a {\em cylinder} with respect to
$M^{l(\overline{y})}$ and generated by the set of solutions
$\varphi(\mathcal{M})$. In any case generating sets for cylinders
coincide their {\em projections}, i.e., sets of solutions for
formulas $\exists\overline{z}\varphi(\overline{x})$, where
$\overline{z}\subset\overline{x}$.}
\end{note}

\begin{proposition}\label{pr031}
If $\mathcal{M}$ is a $n$-element structure, for $n\in\omega$,
then ${\rm ar}({\rm Th}(\mathcal{M}))\leq n$.
\end{proposition}

\begin{proposition}\label{propunfin} {\rm \cite{aafot}.} Any theory of a finite
structure $\mathcal{M}$ is unary-tizable, by an expansion using
finitely many new binary language symbols.
\end{proposition}

\begin{proposition}\label{propunfinfor} {\rm \cite{aafot}.} Any formula of a theory having finitely many solutions is unary-tizable.
\end{proposition}

{\bf Definition} \cite{ananat}.  A theory $T$ is called {\em
almost $n$-ary} if there are finitely many formulae
$\varphi_1(\overline{x}),\ldots,\varphi_m(\overline{x})$ such that
each $T$-formula is $T$-equivalent to a Boolean combination of
$n$-formulae and formulae obtained by substitutions of free
variables in
$\varphi_1(\overline{x}),\ldots,\varphi_m(\overline{x})$.

In such a case we say that the formulae
$\varphi_1(\overline{x}),\ldots,\varphi_m(\overline{x})$ witness
that $T$ is almost $n$-ary.

Almost $1$-ary theories are called {\em almost unary}, almost
$2$-ary theories are called {\em almost binary}, almost $3$-ary
theories are called {\em almost ternary}, etc.

A theory $T$ is called {\em almost $n$-aritizable} if some
expansion $T'$ of $T$ is almost $n$-ary.

Almost $1$-aritizable theories are called {\em almost
unary-tizable}, almost $2$-aritizable theories are called {\em
almost binarizable}, almost $3$-aritizable theories are called
{\em almost ternarizable}, etc.

Assuming that the witnessing set
$\{\varphi_1(\overline{x}),\ldots,\varphi_m(\overline{x})\}$ is
minimal for the almost $n$-ary theory $T$ we have either $m=0$ of
$l(\overline{x})>n$.

Thus we have two minimal characteristics witnessing the almost
$n$-arity of $T$: $m$ and $l(\overline{x})$. The pair
$(m,l(\overline{x}))$ is called the {\em degree} of the almost
$n$-arity of $T$, or the {\em {\rm aar}-degree} of $T$, denoted by
${\rm deg}_{\rm aar}(T)$. Here we assume that $n$ is minimal with
almost $n$-arity of $T$, this $n$ is denoted by ${\rm aar}(T)$.
Clearly, ${\rm aar}(T)\leq{\rm ar}(T)$, and if $m=0$, i.e.,
$n={\rm ar}(T)={\rm aar}(T)$ then it is supposed that
$l(\overline{x})=0$, too.

We have ${\rm aar}(T)\in\omega$ if and only if ${\rm
ar}(T)\in\omega$. So if ${\rm ar}(T)=\infty$ then it is natural to
put ${\rm aar}(T)=\infty$.

\section{Arities and aritizabilities of group theories}

\begin{proposition}\label{arit}
Let $\mathcal{M}$ be an infinite structure,
$\varphi(x_1,\ldots,x_m)$ be a formula satisfying the following
conditions:

$(1)$ for any elements $a_1,\ldots,a_{m-1}\in M$ each result
$\varphi_i(x_i)$ of substitution of these elements to the formula
$\varphi(x_1,\ldots,x_m)$ instead of $m-1$ free variables has
positively many and finitely many solutions in $\mathcal{M}$;

$(2)$ for any $i\leq m$ the formula
$$\psi_i(x_1,\ldots,x_{i-1},x_{i+1},\ldots,x_m):=\exists
x_i\varphi(x_1,\ldots,x_m)$$ has cofinitely many solutions in
$\mathcal{M}$.

Then $\varphi(x_1,\ldots,x_m)$ is not $n$-aritizable for any
$n<m$.
\end{proposition}

Proof. Let $\overline{x}=\langle x_1,\ldots,x_m\rangle$. Assume on
contrary that $\varphi(\overline{x})$ is $n$-arizable for some
$n<m$. Thus $\varphi(\overline{x})$ can be represented, for some
expansion $T'$ of $T={\rm Th}(\mathcal{M})$, by a positive Boolean
combination $\chi(\overline{x})$ of $n$-formulae, written in a
disjunctive normal form
$\bigvee\limits_{j}\bigwedge\limits_{k}\theta_{j,k}(\overline{x}_{jk})$,
$l(\overline{x}_{jk})=n$. For each $j$ we denote by
$\chi_j(\overline{x})$ the formula
$\bigwedge\limits_{k}\theta_{j,k}(\overline{x}_{jk})$. Thus the
set $Z$ of solutions for $\varphi(x_1,\ldots,x_m)$ in the
expansion $\mathcal{M}'$ of $\mathcal{M}$ is represented by
Cartesian and mixed sums and products of the sets of solutions for
$\theta_{j,k}(\overline{x}_{jk})$, as well as by unions $Z_j$ of
the sets of solutions for $\chi_j(\overline{x})$. Now, as in (1),
we substitute elements $a_1,\ldots,a_{m-1}\in M'$ into the
formulae $\varphi(\overline{x})$ and $\chi_j(\overline{x})$. Since
$Z=\bigcup\limits_j Z_j$ and the formulae $\varphi_i(x_i)$ have
finitely many solutions, then the correspondent results
$\chi_{ij}(x_i)$ of substitutions for $\chi_j(\overline{x})$ have
finitely many solutions, too. Now by (1) and (2) we choose $j$
with positively many solutions for $\chi_{ij}(x_i)$ and infinitely
many solutions for conjunctive members
$\theta_{j,k}(\overline{x}_{jk})$, producing the mixed product
$Z_j$ such that varying the tuple $\langle
a_1,\ldots,a_{m-1}\rangle$ copies of $\chi_{ij}(x_i)$ run an
infinite set of solutions. Now we step-by-step show that
intersections of sets of solutions for the formulae
$\theta_{i,j,k}(x_i)$, the results of substitutions of
$a_1,\ldots,a_{m-1}$ into $\theta_{j,k}(\overline{x}_{jk})$, are
infinite, using $l(\overline{x}_{jk})<l(\overline{x})$. It
contradicts the condition that $\varphi_i(x_i)$ has finitely many
solutions in $\mathcal{M}$, as required.

\begin{note} {\rm Notice that the conditions:

$(1)$ a formula $\varphi(x_1,\ldots,x_m)$ has co-infinitely many
solutions in a structure $\mathcal{M}$,

$(2)$ for any $i\leq m$ the formula
$$\psi_i(x_1,\ldots,x_{i-1},x_{i+1},\ldots,x_m):=\exists
x_i\varphi(x_1,\ldots,x_m)$$ has cofinitely many solutions in
$\mathcal{M}$,

\noindent can be realized by unary-tizable, even unary theories.
Indeed, taking a structure $\mathcal{M}$ with a unary predicate
$Q(x)$ with infinitely and co-finitely many solutions, say
$Q=M\setminus\{a\}$ for some $a\in M$, the formula $\varphi(x,y):=
Q(x)\wedge y\approx y$ has co-infinitely many solutions in
$\mathcal{M}$, and both $\exists x \varphi(x,y)$ and $\exists y
\varphi(x,y)$ have cofinitely many solutions in $\mathcal{M}$.}
\end{note}

\begin{theorem}\label{argroup} Let $\mathcal{G}$ be a $\emptyset$-definable subgroup in
a structure $\mathcal{M}$. Then the following conditions are
equivalent:

$(1)$ all formulae of ${\rm Th}(\mathcal{M})$ defining subsets of
finite Cartesian powers of $\mathcal{G}$ are $n$-aritizable for
some fixed natural $n$, and produce $n$-aritizable ${\rm
Th}(\mathcal{G})$;

$(2)$ all formulae of ${\rm Th}(\mathcal{M})$ defining subsets of
finite Cartesian powers of $\mathcal{G}$ are unary-tizable, and
produce unary-tizable ${\rm Th}(\mathcal{G})$;

$(3)$ $\mathcal{G}$ is a finite group.
\end{theorem}

Proof. $(1)\Rightarrow(3)$. Let all formulae of ${\rm
Th}(\mathcal{M})$ defining subsets of Cartesian products of $G$
are $n$-aritizable for some fixed natural $n$ and $\mathcal{G}$ be
an infinite group. We consider the formula
$\varphi(x_1,\ldots,x_n,y):=y\approx x_1\cdot x_2\cdot\ldots\cdot
x_n$. It satisfies the conditions of Proposition \ref{arit} for
the definable substructure $\mathcal{G}$. Indeed, substituting $n$
parameters $a_1,\ldots,a_n$ instead of any $n$ variables we obtain
some unique solution for the rest variable: for
$\varphi(a_1,\ldots,a_n,y)$ we have the solution $y=a_1\cdot
a_2\cdot\ldots\cdot a_n$, and for
$\varphi(a_1,\ldots,a_{i-1},x_i,a_{i+1},\ldots,a_n)$,
$x_i=a^{-1}_{i-1}\cdot\ldots\cdot a^{-1}_{1}a_n
a^{-1}_{n-1}\cdot\ldots\cdot a^{-1}_{i+1}$. Besides, each formula
$\exists x_i\varphi(x_1,\ldots,x_n,y)$, $\exists
y\varphi(x_1,\ldots,x_n,y)$ produces the set $G^n$ of solutions.
By Proposition \ref{arit} the formula $\varphi(x_1,\ldots,x_n,y)$
is not $n$-aritizable. Since $n$ is chosen arbitrary we obtain a
contradiction with the assumption (1).

$(3)\Rightarrow(2)$ follows by Propositions \ref{propunfin} and
\ref{propunfinfor}.

$(2)\Rightarrow(1)$ holds by the monotony of aritizability.

The theorem is proved.


\begin{theorem}\label{coargroup} For any group $\mathcal{G}$ the following conditions
are equivalent:

$(1)$ ${\rm Th}(\mathcal{G})$ is aritizable,

$(2)$ ${\rm Th}(\mathcal{G})$ is almost $n$-aritizable for some
$n$,

$(3)$ ${\rm Th}(\mathcal{G})$ is unary-tizable,

$(4)$ ${\rm Th}(\mathcal{G})$ is almost unary-tizable,

$(5)$ ${\rm Th}(\mathcal{G})$ is $n$-ary for some $n\in\omega$,

$(6)$ $\mathcal{G}$ is a finite group.
\end{theorem}

Proof. $(1)\Leftrightarrow(2)$ holds, since if a theory $T$ is
aritizable then it is $n$-aritizable for some $n$ and therefore
almost $n$-aritizable. Conversely, if a theory $T$ is almost
$n$-aritizable for some $n$ which is witnessed by a finite set
$\Phi$ of $T$-formulae, then $T$ is $m$-aritizable, where $m={\rm
max}\{n,l(\Phi)\}$, $l(\Phi)$ is a maximal length
$l(\overline{x})$ of free variable tuple $\overline{x}$, for
$\varphi(\overline{x})\in\Phi$.

$(1)\Leftrightarrow(3)\Leftrightarrow(6)$ follows by Theorem
\ref{argroup}.

$(6)\Rightarrow(5)$ holds by Proposition \ref{pr031}.

$(3)\Rightarrow(4)$, $(4)\Rightarrow(2)$, and $(5)\Rightarrow(1)$
are obvious.

The theorem is proved.

\medskip
Theorem \ref{coargroup} immediately implies the following
dichotomy:

\begin{corollary}\label{coarfininfgr}
For any group $\mathcal{G}$ either ${\rm ar}({\rm
Th}(\mathcal{G}))\in\omega$, if $\mathcal{G}$ is finite, or ${\rm
ar}({\rm Th}(\mathcal{G}))=\infty$, if $\mathcal{G}$ is infinite.
\end{corollary}

\section{Aritizable theories of monoids and groupoids}

\medskip
Arguments for Theorem \ref{argroup} and for Corollary
\ref{coargroup} confirm that monoid and groupoid theories
correspondent to a theory of a group $G$ are non-aritizable if $G$
is infinite. At the same time the following example shows that
Corollaries \ref{coargroup} and \ref{coarfininfgr} do not hold for
monoids.

\begin{example}\label{exarinf} {\rm Let $\mathcal{M}$ be an infinite monoid with the unit $e$ and the
rule
$$\forall x,y(x\ne e\wedge y\ne e\to x\cdot y = e).$$ The formula
$\varphi(x,y,z):=(x\cdot y=z)$, defining the multiplication in
$\mathcal{M}$, is ${\rm Th}(\mathcal{M})$-equivalent to the
following Boolean combination of binary formulae:
$$
(x=e\wedge y=z)\vee(y=e\wedge x=z)\vee(x\ne e\wedge y\ne e\wedge
z=e),
$$
witnessing that ${\rm Th}(\mathcal{M})$ is unary and implying that
${\rm Th}(\mathcal{M})$ is unary-tizable.}
\end{example}

\begin{note}\label{notarmon}
{\rm Similarly Example \ref{exarinf}, Corollaries \ref{coargroup}
and \ref{coarfininfgr} fail for monoids $\mathcal{M}$ with finite
$R(\mathcal{M})=\{a\cdot b\mid a,b\in M, a\ne e, b\ne e \}$, and
for related groupoids.

Indeed, since $R(\mathcal{M})$ is finite consisting of elements
$c_1,\ldots,c_m$, there is an expansion $\mathcal{M}'$ of
$\mathcal{M}$ by finitely many disjoint binary predicates
$D_i=\{\langle a,b\rangle\mid a\ne e, b\ne e, a\cdot b=c_i\}$, and
unary single predicates $R_i=\{c_i\}$, $i\leq m$, with
$D_1\cup\ldots\cup D_m=(M\setminus\{a\})^2$. The formula
$\varphi(x,y,z):=(x\cdot y=z)$, defining the multiplication in
$\mathcal{M}$, is ${\rm Th}(\mathcal{M})$-equivalent to the
following Boolean combination of binary formulae:
$$
(x=e\wedge y=z)\vee(y=e\wedge x=z)\vee\left(x\ne e\wedge y\ne
e\wedge\left(\bigvee\limits_{i=1}^m(D_i(x,y)\wedge
R_i(z))\right)\right).
$$
witnessing \ that \ ${\rm Th}(\mathcal{M}')$ \ is \ binary \ and \
implying \ that \ ${\rm Th}(\mathcal{M})$ \ is \ binarizable, \
with \ ${\rm aar}({\rm Th}(\mathcal{M}))=1$.

The same effect of binarizability is satisfied for groupoids
$\mathcal{M}$ with finite $R(\mathcal{M}):=M\cdot M$.}
\end{note}

In view of Remark \ref{notarmon} the following assertion holds:

\begin{proposition}\label{prmongrbin} \ Any \ monoid \ {\rm (}groupoid{\rm
)} \ $\mathcal{M}$ \ with \ finite \ $R(\mathcal{M})$ \ is \
binarizable \ with \ ${\rm aar}({\rm Th}(\mathcal{M}))=1$.
\end{proposition}

Arguments for Remark \ref{notarmon} show that Proposition
\ref{prmongrbin} can be spread for each algebra
$\mathcal{A}=\langle A;f^{n_1}_1,\ldots,f^{n_k}_k\rangle$ with
finite $R(\mathcal{A})=\bigcup\limits_{i=1}^k f_i(A,A,\ldots,A)$.
In such a case we replace binary relations $D_i$ by $n_j$-ary
relations $D^{n_j}_{i,j}=\{\langle a_1,\ldots,a_{n_j}\rangle\mid
f_j(a_1,\ldots,a_{n_j})=c_i\}$, $1\leq j\leq k$.

\begin{proposition}\label{pralg} Any algebra $\mathcal{A}$ in a finite language and with finite  $R(\mathcal{A})$ is
aritizable and ${\rm aar}({\rm Th}(\mathcal{A}))=1$.
\end{proposition}

\begin{note}
{\rm Since arities and aritizabilities are preserved under
disjoint unions and compositions of theories \cite{aafot, ananat}
Proposition \ref{pralg} admits natural generalizations for
theories in infinite languages.}
\end{note}

\section{Conclusion}

We considered possibilities for arities of group, monoid and
groupoid theories. It is shown that both $n$-ary and aritizable
group theories are exactly theories of finite groups. This
property fails for theories of monoids, theories of groupoids, and
theories of some universal algebras. A natural question arises on
a characterization of arities and aritizabilities of theories of
various algebras.

\noindent Novosibirsk State Technical University, \\ 20, K.Marx
avenue, Novosibirsk, 630073, Russia; \\ Sobolev Institute of
Mathematics, \\ 4, Acad. Koptyug avenue, Novosibirsk, 630090,
Russia; \\ Novosibirsk State University, \\ 1, Pirogova street,
Novosibirsk, 630090, Russia

\medskip\noindent
e-mail: inessa7772@mail.ru, sudoplat@math.nsc.ru

\end{document}